\documentclass[12pt]{amsart}
\usepackage{latexsym}

\usepackage{amsfonts}
\usepackage{amsthm}
\usepackage{amsmath}
\usepackage{amssymb}

\def\sw#1{{\sb{(#1)}}}
\def\<{{\langle}}
\def\>{{\rangle}}

\def\eps{\varepsilon}

\def\note#1{{}}

\def\note#1{}

\def\beq{\begin{equation}}
\def\eeq{\end{equation}}

\def\id{{I}}

\def\ot{{\otimes}}

\newcounter{zlist}

\def\Label#1{\label{#1}\ifmmode\llap{[#1] }\else
\marginpar{\smash{\hbox{\tiny [#1]}}}\fi}
\def\Label{\label}

\newtheorem{proposition}{Proposition}[section]

\newtheorem{theorem}[proposition]{Theorem}

\theoremstyle{definition}
\newtheorem{definition}[proposition]{Definition}

\theoremstyle{remark}
\newtheorem{remark}[proposition]{Remark}

\newcounter{c}

\newcommand{\etyk}[1]{\vspace{-7.4mm}$$\begin{equation}\Label{#1}
\addtocounter{c}{1}}
\renewcommand{\]}{\ifnum \value{c}=1 $$\else \end{equation}\fi}
\begin{document}

\title{Algebra structures arising from Yang-Baxter systems}
\author{Barbu R. Berceanu}
\author{Florin F. Nichita}
\author{Calin Popescu}
\address{Abdus Salam School of Mathematical Sciences, GC University, 
Lahore, Pakistan}
\address{Institute of Mathematics of the Rumanian Academy, 
P.O. Box 1-764, RO-70700, Bucharest, Rumania}
\email{Barbu.Berceanu@imar.ro}
\address{Institute of Mathematics of the Rumanian Academy, 
P.O. Box 1-764, RO-70700, Bucharest, Rumania}
\email{Florin.Nichita@imar.ro}
\email{Calin.Popescu@imar.ro}
\date{May, 2010}
\begin{abstract}
Yang-Baxter operators from algebra structures
appeared for the first time in \cite{Nic:sel}, \cite{nu} and \cite{DasNic:yan}.
Later, Yang-Baxter systems from entwining structures were constructed
in \cite{BrzNic:coa}.
In this paper we show
that an algebra factorisation 
can be constructed from
a Yang-Baxter system.
\end{abstract}
\maketitle

\section{Introduction}

The quantum Yang-Baxter equation plays a crucial role in
theoretical physics, 
in the theory of quantum groups,
in knot theory,  in  the theory of braided monoidal categories, etc. 
Due to the complexity of various integrable models the need arose for extensions and generalisations of Yang-Baxter equations and related algebraic structures.  One such  generalisation, termed a {\em Yang-Baxter system}, was proposed in  \cite{HlaSno:sol} in the context of non-ultralocal integrable systems  \cite{HlaKun:qua}. 

Motivated by the need for developing a general theory of non-commutative principal bundles, entwining structures were introduced in \cite{BrzMaj:coa} as generalised symmetries of such bundles. In recent years, entwining structures and categories of modules associated to them have been used to provide a unification of various types of Hopf modules (cf.\ \cite{Brz:mod}, \cite{CaeMil:Fro}), and eventually have led to the revival of the theory of corings (cf.\ \cite{Brz:str}, \cite{BrzWis:cor}).

Thus, Yang-Baxter systems and entwining structures arose in entirely different contexts, were motivated by as far fields as high energy physics and non-commutative geometry, and have entirely different applications. There is no reason to expect that there is any connection between these two notions. 
Yet, in \cite{BrzNic:coa}, we proved that such a connection, and a very close one, indeed exists. More precisely, we showed that to any entwining
structure one can associate a Yang-Baxter system.

In this paper, we consider another construction. 
We start
with a Yang-Baxter system and construct an algebra factorisation.
As observed in \cite[Proposition~2.7]{BrzMaj:coa}, entwining structures are related to {\em algebra factorisations} via semi-dualisation (see \cite[p.\ 300]{Maj:fou} for discussion of algebra factorisations). Note that
Yang-Baxter operators from algebra structures
appeared for the first time in \cite{Nic:sel}, \cite{nu} and \cite{DasNic:yan}.
Other interesting papers
related to the constructions appearing in this article are: \cite{AshBer:sim}, \cite{GraOgi:bra}, etc.

The next section contains preliminaries on the Yang-Baxter equation. 
The main results  are contained in Section 3. 

\section{Preliminaries and Review}

Throughout this paper $k$ is a field. 

The unadorned tensor product is over $k$. 

The identity map on a $k$-vector space $V$ is denoted by $\id_V$ or simply by $\id$ provided the domain is clear from the context. All algebras are over $k$, 
they are associative and with unit 1.  
The product in an algebra $A$ is denoted by $\mu:A\ot A\to A$, while the unit map is denoted by $\iota:k\to A$. 

All coalgebras are over $k$, they are coassociative  and with a counit. 
 Coproduct in a coalgebra $C$ is denoted by $\Delta:C\to C\ot C$ and the counit by  $\eps:C\to k$.   %We assume that ${\rm dim}\, C \geq 1$.  
We use the standard notation for coalgebras. In particular, for a coalgebra $C$, we use Sweedler's notation to denote the coproduct $\Delta$ on elements, i.e., $\Delta(c) = \sum c\sw1\ot c\sw 2$, for all $c\in C$.

\bigskip

For any vector spaces $ V $ and $ W $,
$ \tau_{V, W}: V \otimes W \rightarrow W \ot V \  $ denotes the natural
bijection defined by $ \tau_{V, W}(v \ot w) = w \ot v $.

Let $ R: V \ot V \rightarrow V \ot V  $
be a $ k$-linear map. Define
$ {R_{12}}= R \ot \id_V$, ${R_{23}}= \id_V \ot R$  and
${R_{13}}=(\id_V\ot\tau_{V, V})\circ(R\ot \id_V)\circ (\id_V\ot \tau_{V, V})$. Each of the $R_{ij}$ is thus a linear endomorphism of $V\ot V\ot V$.

\begin{definition}\label{def.flo}
 An invertible  $k$-linear map  $ R : V \ot V \rightarrow V \ot V $
is called a {\em Yang-Baxter
operator} (or simply a {\em YB operator}) if it satisfies the  equation
\begin{equation}  \label{ybeq}
R_{12}  \circ  R_{23}  \circ  R_{12} = R_{23}  \circ  R_{12}  \circ  R_{23}.
\end{equation}

\end{definition}

Equation (\ref{ybeq}) is usually called the {\em braid equation}. It is a
well-known fact that the operator $R$ satisfies (\ref{ybeq}) if and only if
$R\circ \tau_{V, V}  $ satisfies
   the {\em quantum Yang-Baxter equation} 
\begin{equation}   \label{ybeq2}
R_{12}  \circ  R_{13}  \circ  R_{23} = R_{23}  \circ  R_{13}  \circ  R_{12}.
\end{equation}
For a review of Yang-Baxter operators we refer to \cite{lara}, \cite{Maj:fou}
and \cite{Nic:non}.

\bigskip

In what follows, we use the following construction of Yang-Baxter operators described in \cite{DasNic:yan}.

If $A$ is a $k$-algebra, then for all non-zero $r,s\in k$, the linear map 
\begin{equation}\label{ra}
R^A_{r,s}:A \ot A
\rightarrow A \ot A,\quad   a \ot b \mapsto s ab \ot 1 + r 1 \ot ab - s a \ot b 
\end{equation}
is a Yang-Baxter operator. 

The inverse of $R^A_{r,s}$ is 
${(R^A_{r,s})}^{-1}(a \ot b)= \frac{1}{r} ab \ot 1 + \frac{1}{s} 1 \ot 
ab - \frac{1}{s} a \ot b $.

\bigskip

Dually, if $C$ is a coalgebra, then for all non-zero $p,t\in k$, the linear map 
\begin{equation}\label{rc}
R_C^{p,t}:C \ot C
\rightarrow C \ot C,\quad   c \ot d \mapsto p \eps (c) \sum d\sw 1 \ot d\sw 2 + t \eps (d) 
\sum c\sw 1 \ot c\sw 2 - p c \ot d 
\end{equation}
is a Yang-Baxter operator. 

\bigskip

Note that in both cases the assumption that parameters are non-zero is needed only for the invertibility of $R^A_{r,s}$ and $R_C^{p,t}$. $R^A_{r,s}$ and $R_C^{p,t}$ satisfy the braid relation for any value of $r,s,p$ and $t$.

\bigskip

%\section{From entwining structures to Yang-Baxter systems (and back)}

Yang-Baxter systems were introduced in \cite{HlaSno:sol} as a spectral-parameter independent generalisation of quantum Yang-Baxter equations related to non-ultralocal integrable systems studied previously in  \cite{HlaKun:qua}. 

Yang-Baxter systems are conveniently defined in terms of {\em Yang-Baxter commutators}. Consider three vector spaces $V,V',V''$ and  three linear maps
$ R : V \ot V' \rightarrow V \ot V' $,
$ S : V \ot V'' \rightarrow V \ot V'' $ and
$ T : V' \ot V'' \rightarrow V' \ot V'' $. Then a {\em Yang-Baxter commutator} is a map
$ [R,S,T]:  V \ot V' \ot V'' \rightarrow V \ot V' \ot V'' $, defined by
\begin{equation}   \label{ybcomm}
[R,S,T]= R_{12}  \circ  S_{13}  \circ  T_{23} - T_{23}  \circ  S_{13} 
\circ  R_{12} \ .
\end{equation}
In terms of a Yang-Baxter commutator,  the quantum Yang-Baxter equation (\ref{ybeq2})  is expressed simply as $[R,R,R] = 0$.

\begin{definition}\label{def.wxz}
Let $V$ and  $V'$ be vector spaces. A system of linear maps 
$$ W : V \ot V \rightarrow V \ot V ,\quad 
Z : V' \ot V' \rightarrow V' \ot V' , \quad
 X : V \ot V' \rightarrow V \ot V' 
$$ 
is called  a {\em WXZ-system} or a 
{\em Yang-Baxter system}, provided the following equations are satisfied:
\begin{equation}   \label{ybeqn4}
[W,W,W]\ = \  0 \ ,
\end{equation}
\begin{equation}   \label{ybeqn5}
[Z,Z,Z]\ = \  0 \ ,
\end{equation}
\begin{equation}   \label{ybeqn6}
[W,X,X]\ = \  0 \ ,
\end{equation}
\begin{equation}   \label{ybeqn7}
[X,X,Z]\ = \  0 \ .
\end{equation}
\end{definition}
\bigskip

There are several algebraic origins and applications of WXZ-systems (see, for example, \cite{Vla:met}, \cite{Sno:con}, etc).

\begin{remark}
Given a WXZ-system as in Definition~\ref{def.wxz} one can construct a Yang-Baxter operator on $V\oplus V'$, provided the map $X$ is invertible. This is a special case of a {\em gluing procedure} described in  \cite[Theorem~2.7]{MajMar:glu} (cf.\  \cite[Example 2.11]{MajMar:glu}).
Let $ R = W \circ \tau_{V,V} $, $ R' = Z \circ \tau_{V',V'} $, $ U = X
\circ \tau_{V',V} $. Then the linear map
$$
 R \oplus_{U} R': (V \oplus V') \ot (V \oplus V') \to  (V \oplus V') \ot (V \oplus V')
$$ 
given by
$ R \oplus_{U} R'|_{V \ot V} = R $,
$ R \oplus_{U} R'|_{V' \ot V'} = R' $, and for all $ x \in V $, $ y \in V' $,
$$ 
(R \oplus_{U} R')(y \ot x) = U(y \ot x ), \qquad  (R \oplus_{U} R')(x \ot y ) = U^{-1}(x \ot y)
$$
is a Yang-Baxter operator.
\end{remark}

\bigskip

Entwining structures were introduced in \cite{BrzMaj:coa} in order to recapture the symmetry structure of non-commutative (coalgebra) principal bundles or coalgebra-Galois extensions. For  applications and algebraic content of entwining structures we refer to  \cite{CaeMil:Fro} and \cite{BrzWis:cor}.

\begin{definition}\label{def.entw} 
 An algebra $A$ is said to be {\em entwined} with  a coalgebra $C$ if there exists  a linear  map 
$ \psi : C \ot A
\rightarrow A \ot C $ satisfying the following four conditions:

(1) $ \psi \circ (I_{C} \ot \mu) = ( \mu \ot I_{C}) \circ (I_{A} \ot
\psi) \circ ( \psi \ot I_{A}) $,

(2) $ (I_{A} \ot \Delta) \circ \psi = ( \psi \ot I_{C} ) \circ (I_{C}
\ot \psi) \circ ( \Delta \ot I_{A} ) $,

(3) $ \psi \circ ( I_{C} \ot \iota ) = \iota \ot I_{C} $,

(4) $ (I_{A} \ot \eps ) \circ \psi = \eps \ot I_{A} $.

The map $ \psi $ is known as an {\em entwining map}, and the triple  $ {(A,C)}_{\psi} $ is called an {\em entwining structure}. 
\end{definition}

To denote the action of an entwining
map $ \psi $ on elements it is convenient to use the following {\em $\alpha$-notation}, for all $a,b\in A$ and $c\in C$,
$$
\psi (c \ot a) = \sum_{\alpha} a_{\alpha} \ot c^{\alpha} , \quad
(I_{A} \ot \psi) \circ ( \psi \ot I_{A}) (c \ot a \ot b) =
\sum_{\alpha, \beta} a_{\alpha} \ot b_{\beta} \ot c^{\alpha \beta},
$$
etc. The relations (1), (2), (3)
and (4) in Definition~\ref{def.entw} are  equivalent to the following explicit relations, for all $ a,b
\in A$, $c \in C$,
\begin{equation}\label{a}
\sum_{\alpha}(ab)_{\alpha}\ot c^{\alpha} = \sum_{\alpha, \beta}
a_{\alpha} b_{\beta} \ot c^{\alpha \beta},
\end{equation}
\begin{equation}\label{b}
\sum_{\alpha} a_{\alpha} \ot {c^{\alpha}}\sw{1} \ot {c^{\alpha}}\sw{2} =
\sum_{\alpha, \beta} a_{\beta \alpha} \ot {c\sw{1}}^{\alpha} \ot
{c\sw{2}}^{\beta},
\end{equation}
\begin{equation}\label{c}
\sum_{\alpha} 1_{\alpha} \ot c^{\alpha} = 1 \ot c,
\end{equation}
\begin{equation}\label{d}
\sum_{\alpha}a_{\alpha} \eps (c^{\alpha}) = a \eps (c).
\end{equation}

\begin{theorem}\label{thm.main} (\cite{BrzNic:coa})
 Let $A$ be an algebra  and let $C$ be a
coalgebra. For any $ s, r, t, p \in k$ define linear maps 
$$W: A\ot A\to A\ot A, \qquad a\ot b \mapsto s ba \ot 1 + r 1 \ot ba - s b \ot a, 
$$
$$
Z: C\ot C\to C\ot C, \qquad c\ot d \mapsto t \eps (c) \sum d\sw 1 \ot d\sw 2 + p \eps (d) 
\sum c\sw 1 \ot c\sw 2 - p d \ot c.
$$
Let $X:A\ot C\to A\ot C$ be a linear map such that $ X\circ (\iota \ot \id_C) = \iota\ot \id_{C} $ and $ (\id_A \ot \eps ) \circ X =\id_A\ot \eps$. Then $W,X,Z$ is a Yang-Baxter system if and only if $A$ is entwined with $C$ by the map $\psi: = X\circ \tau_{C, A}$.

\end{theorem}

\begin{remark}
As observed in \cite[Proposition~2.7]{BrzMaj:coa}, entwining structures are related to {\em algebra factorisations} via semi-dualisation (see \cite[p.\ 300]{Maj:fou} for discussion of algebra factorisations). The arguments similar to those in the proof of Theorem~\ref{thm.main} show that, given two algebras $A$, $B$, and an algebra factorisation map $\Psi: B\ot A\to A\ot B$, one can construct a Yang-Baxter system with $W$ and $Z$ of the same form as $W$ in Theorem~\ref{thm.main} and $X=\Psi\circ\tau_{A,B}$. 
\end{remark}

\section{Main Result and Consequences}

We give the formal definition of 
 {\em algebra factorisations}
below.

\begin{definition}\label{def.algfact} 
 An algebra $A$ is said to be {\em entwined} with  an algebra $B$ if there exists  a linear  map 
$ \phi : B \ot A
\rightarrow A \ot B $ satisfying the following four conditions:

(1) $ \phi \circ (I_{B} \ot \mu_A ) = ( \mu_A \ot I_{B}) \circ (I_{A} \ot
\phi) \circ ( \phi \ot I_{A}) $,

(2) $ \phi \circ (\mu_B \ot I_A) = ( I_{A} \ot \mu_B) \circ (
\phi \ot I_B) \circ ( I_B \ot \phi) $,

(3) $ \phi \circ ( I_{B} \ot \iota_A ) = \iota_A \ot I_{B} $,

(4) $ \phi \circ ( \iota_B \ot I_{A} ) = I_{A} \ot \iota_B $,

The map $ \phi $ is known as an {\em algebra factorisation map}, and the triple  $ {(A,B)}_{\phi} $ is called an {\em algebra factorisation}. 
\end{definition}

\bigskip

Let 
$R:V\otimes V\rightarrow V\otimes V$ 
be a Yang-Baxter operator (see Definition~\ref{def.flo})
and
let the operators   
$\sigma_i = R_{i,i+1}: V^{\ot n+1} \rightarrow  V^{\ot n+1} $, for
$i=1,2,3,\cdots n$.
 
\bigskip

We define a {\it braid product} on the tensor 
algebra 
$T(V)$ 
as follows: 

for 
$x = x_1\otimes\cdots\otimes x_n\in 
T^n(V)$, 
$n\geq 1$, 
and 
$y = y_1\otimes\cdots\otimes y_p\in 
T^p(V)$, 
$p\geq 1$, 
set 

$$
xy = 
(\sigma_p\cdots\sigma_{n+p-1})\cdots 
(\sigma_1\cdots\sigma_n)
(x\otimes y); 
$$

that is, 
$\mu_{n,p}:T^n(V)\otimes T^p(V)
\rightarrow T^{n+p}(V)$, 
$$
\mu_{n,p} = \prod_{i=0}^{p-1}
(\sigma_{p-i}\cdots\sigma_{n+p-i-1}).
$$

\begin{theorem}
Starting with a  WXZ-system as in Definition~\ref{def.wxz}, and
$ R = W \circ \tau_{V,V} $, $ R' = Z \circ \tau_{V',V'} $, $ U = X
\circ \tau_{V',V} $, we have the following properties.

a) $ ( \  T(V), \  \mu, \  1_k \ ) $ and
 $ ( \ T(V'), \  \mu', \  1_k \ ) $ are k-algebra structures.

b) There exists an algebras factorisation 
$ \psi : T(V') \otimes T(V) \rightarrow T(V) \otimes T(V')$, 
defined by

$ \psi_{n,m} : T^n(V') \otimes T^m(V) \rightarrow T^m(V) \otimes T^n(V'),$

$ x \otimes y \mapsto (U_{m, m+1} \cdots U_{n+m-1, m+n})\cdots 
( U_{12}\cdots U_{n,n+1})
(x\otimes y) $

$ x \otimes 1_k \mapsto 1_k \otimes x $

$ 1_k \otimes y \mapsto y \otimes 1_k $

\end{theorem}

{ \bf Proof.}

We prove the first part of the theorem, part a),
 below.

The operators   
$\sigma_i = R_{i,i+1}$, 
$i=1,2,3,\cdots$, 
have the following properties: 
\begin{itemize}
\item[{\bf (a)}] 
$\sigma_i\sigma_{i+1}\sigma_i = 
\sigma_{i+1}\sigma_i\sigma_{i+1}$, 
for $i=1,2,3,\cdots$; 
\item[{\bf (b)}] 
$\sigma_i\sigma_j = \sigma_j\sigma_i$, 
whenever 
$|i - j| > 1$;
\item[{\bf (c)}] 
if 
$i\leq k < j$, 
then 
$(\sigma_i\cdots\sigma_j)\sigma_k = 
\sigma_{k+1}(\sigma_i\cdots\sigma_j)$ --- 
this follows from {\bf (a)}; 
and  
\item[{\bf (d)}] 
if 
$i < j$, 
and  
$w(\sigma_i,\cdots,\sigma_{j-1})$
is any word in 
$\sigma_i$, $\cdots$, $\sigma_{j-1}$, 
then 
$$
(\sigma_i\cdots\sigma_j)
w(\sigma_i,\cdots,\sigma_{j-1}) = 
w(\sigma_{i+1},\cdots,\sigma_j)
(\sigma_i\cdots\sigma_j). 
$$
(This follows from {\bf (c)}.)
\end{itemize} 

The {\it braid product} on the tensor 
algebra 
$T(V)$ 
is defined as follows: For 
$x = x_1\otimes\cdots\otimes x_n\in 
T^n(V)$, 
$n\geq 1$, 
and 
$y = y_1\otimes\cdots\otimes y_p\in 
T^p(V)$, 
$p\geq 1$, 
set 
$$
xy = 
(\sigma_p\cdots\sigma_{n+p-1})\cdots 
(\sigma_1\cdots\sigma_n)
(x\otimes y); 
$$
that is, 
$\mu_{n,p}:T^n(V)\otimes T^p(V)
\rightarrow T^{n+p}(V)$, 
$$
\mu_{n,p} = \prod_{i=0}^{p-1}
(\sigma_{p-i}\cdots\sigma_{n+p-i-1}).
$$ 

Associativity, 
$\mu_{n+p,q}\mu_{n,p} = 
\mu_{n,p+q}\mu_{p,q}$, 
amounts to  
$$
\begin{aligned}
\left(
\prod_{i=0}^{q-1}
(\sigma_{q-i}\cdots\sigma_{n+p+q-i-1})
\right)  
\left(
\prod_{i=0}^{p-1}
(\sigma_{p-i}\cdots\sigma_{n+p-i-1})
\right) = & \\
\left(
\prod_{i=0}^{p+q-1}
(\sigma_{p+q-i}\cdots\sigma_{n+p+q-i-1})
\right)  
\left(
\prod_{i=0}^{q-1}
(\sigma_{n+q-i}\cdots\sigma_{n+p+q-i-1})
\right); &    
\end{aligned}
\eqno{(1)}
$$
more explicitly, 
$$
\begin{aligned} 
\underbrace{
(\sigma_q\cdots\sigma_{n+p+q-1})\cdots 
(\sigma_1\cdots\sigma_{n+p})}_{q \,\, 
\textrm{ $(n+p)$-blocks}} 
\underbrace{
(\sigma_p\cdots\sigma_{n+p-1})\cdots 
(\sigma_1\cdots\sigma_n)}_{p \,\, 
\textrm{ $n$-blocks}} = & \\ 
\underbrace{
(\sigma_{p+q}\cdots\sigma_{n+p+q-1})\cdots 
(\sigma_1\cdots\sigma_n)}_{p+q \,\,
\textrm{ $n$-blocks}} 
\underbrace{ 
(\sigma_{n+q}\cdots\sigma_{n+p+q-1})\cdots 
(\sigma_{n+1}\cdots\sigma_{n+p})}_{q \,\,
\textrm{ $p$-blocks}}. &  
\end{aligned}
\eqno{(1')}
$$
To prove $(1)$, we use {\bf (b)} and 
{\bf (d)} repeatedly. First, use {\bf (d)} 
successively for 
$i = 0, \cdots, p - 1$, 
to shift each $n$-block 
$$
\sigma_{p-i}\cdots\sigma_{n+p-i-1}
$$ 
$q$ times leftwards and thus get 
$$
\begin{aligned}
\left(
\prod_{i=0}^{q-1}
(\sigma_{q-i}\cdots\sigma_{n+p+q-i-1})
\right)  
\left(
\prod_{i=0}^{p-1}
(\sigma_{p-i}\cdots\sigma_{n+p-i-1})
\right) = & \\
\left(
\prod_{i=0}^{p-1}
(\sigma_{p+q-i}\cdots\sigma_{n+p+q-i-1})
\right)  
\left(
\prod_{i=0}^{q-1}
(\sigma_{q-i}\cdots\sigma_{n+p+q-i-1})
\right) . &    
\end{aligned}
\eqno{(2)}
$$
Next, transform the second product in the 
right-hand member of $(2)$ as follows: 
Split each $(n + p)$-block into an $n$-block 
followed by a $p$-block: 
$$
\sigma_{q-i}\cdots\sigma_{n+p+q-i-1} = 
(\sigma_{q-i}\cdots\sigma_{n+q-i-1})
(\sigma_{n+q-i}\cdots\sigma_{n+p+q-i-1}).
$$
By {\bf (b)}, for 
$i = 0, \cdots, q - 1$, 
the $i$-th such $n$-block in the product 
commutes with the $i$ such $p$-blocks 
preceding it. The $q$-step procedure 
consisting, at step $i$, in moving the 
$i$-th $n$-block to the left of the now 
contiguous sequence of $i$ $p$-blocks that 
precede it, yields
$$
\begin{aligned}
 & \prod_{i=0}^{q-1}
(\sigma_{q-i}\cdots\sigma_{n+p+q-i-1}) \\ 
 & = \left( 
\prod_{i=0}^{q-1}
(\sigma_{q-i}\cdots\sigma_{n+q-i-1})
\right) 
\left( 
\prod_{i=0}^{q-1}
(\sigma_{n+q-i}\cdots\sigma_{n+p+q-i-1})
\right). 
\end{aligned}    
\eqno{(3)}
$$
Relations $(2)$ and $(3)$ yield the desired 
result. 

The next part of the theorem, part b), follows in a similar manner.
It will be included in the extended version of this paper 
(and follows from \cite{AshBer:sim}).

\qed

\begin{remark} From the previous theorem,
 $ T(V) \otimes T(V')$ has an algebra structure associated with that
algebra factorisation.
\end{remark}

\begin{remark}
 In a natural way,
 $ ( \  T(V), \  \mu, \  1_k \ ) $ gets a Hopf algebra structure.
The maps $ \Delta(v) = v \ot 1_k + 1_k \ot v $, 
$ \ \eps(v) = 0 $ and $ \ S(v)= -v$ $ \  \forall v \in V$ extend
naturally on $ T(V)$. 
\end{remark}

\end{document}